\documentclass[12pt,reqno]{amsart}
\usepackage{amssymb}

\usepackage[all]{xy}
\usepackage{graphicx,amssymb,amsfonts,epsfig,amsthm,a4,amsmath,url}
\usepackage[latin1]{inputenc}

\swapnumbers

\newcounter{ictr}

\newcounter{actr}

\newtheorem{thm}{Theorem}[section]

\newtheorem{thms}{Theorem}[section]

\newtheorem*{corr}{Corollary}

\newtheorem{finite-union}[thm]{Finite Union Theorem}

\theoremstyle{definition}

\newtheorem{defn}[thm]{Definition}

\newtheorem{defns}[thms]{Definition}

\theoremstyle{remark}

\newtheorem{open}[thm]{Open Question}

\numberwithin{equation}{section}

\newcommand{\hhs}{{Hilbert-Hadamard space}}%{{non-positively curved Hilbert manifold}}

\newcommand{\GL}{\text{GL}}

\title{The Novikov conjecture  
}

\dedicatory{Dedicated to Sergei Novikov on the occasion of his 80th birthday}

\author{Guoliang Yu}
\address{Department of Mathematics, Mailstop 3368, Texas A\&M University, College Station, TX 77843, USA}
\email{guoliangyu@math.tamu.edu}

\thanks{The author is partially supported by  a grant from the
U.S. National Science Foundation.}

\begin{document}

\begin{abstract}
 We give a survey on  recent development of the Novikov conjecture
 and its applications to topological rigidity and non-rigidity.

  \end{abstract}

\maketitle

% \tableofcontents

\section{Introduction}

A central problem in mathematics is the Novikov conjecture \cite{N}. Roughly speaking,
the Novikov conjecture claims that compact smooth manifolds are rigid at an infinitesimal level. More precisely, 
the Novikov conjecture states that the higher signatures of compact oriented smooth manifolds are invariant under orientation-preserving homotopy equivalences.  Recall that a compact manifold is called aspherical if its universal cover is contractible.
In the case of aspherical manifolds, the Novikov conjecture is an infinitesimal version of the Borel conjecture, which states that all compact aspherical manifolds
are topologically rigid, i.e. if another compact manifold $N$ is homotopy equivalent to the given compact aspherical manifold $M$, then $N$ is homeomorphic to $M$.  A deep theorem of Novikov says that the rational Pontryagin classes are invariant under orientation-preserving homeomorphisms  
 \cite{N1}. The Novikov conjecture for compact aspherical manifolds follows from the Borel conjecture and Novikov's theorem since for aspherical manfolds, the information about higher signatures is equivalent to that of rational Pontryagin classes.

The Novikov conjecture   has inspired a lot of beautiful mathematics. It motivated the development of  Kasparov's $KK$-theory, Connes' cyclic cohomology theory, Gromov-Connes-Moscovici theory of almost flat bundles, Connes-Higson's $E$-theory, and quantitative operator $K$-theory. The Novikov conjecture has been proven for a large number of cases. The general philosophy is that the conjecture should be true  if the fundamental group of the manifold arises from nature.
The purpose of this article is to give a friendly survey on recent development of the Novikov conjecture and its applications to topological rigidity and non-rigidity.

There are two natural approaches to the Novikov conjecture: the analytic approach based on ideas from noncommutative geometry and the topological approach using $L$-theory. Up until this point, the analytic approach has been more successful partly due to the crucial fact that it is easier to do cutting and pasting in the world of noncommutative geometry. Cutting and pasting makes it possible to retrieve the $K$-theoretic information necessary to prove the Novikov conjecture.
 In noncommutative geometry, we are dealing
with group $C^\ast$-algebras, while in topology, we are studying group rings. To illustrate this fundamental point, let us consider the case when the fundamental group is the group of integers. In this case,
 the group $C^\ast$-algebra is the algebra of continuous functions on the circle while the group ring is  the ring of Laurent polynomials. We can apply cutting and paste to a continuous function to obtain another continuous functions  using partitioning of unity,
 this procedure cannot be applied to  Laurent polynomials.

This survey  focuses on the analytic approach to the Novikov conjecture
and as well as methods inspired by the analytic approach.
The analytic  approach  uses operator $K$-theoretic techniques from noncommutative geometry \cite{M, BC, C, K}. The Novikov conjecture follows from the rational strong Novikov conjecture, which states that the rational Kasparov-Baum-Connes map is injective, here the Kasparov-Baum-Connes map is the assembly map from $K$-homology group of the classifying space for the fundamental group to the $K$-theory of group $C^\ast$-algebra associated to the
fundamental group.
Using this approach, the Novikov conjecture has been proven when the fundamental group of the manifold belongs one of the following cases:
\begin{enumerate}
	\item groups acting properly and isometrically on simply connected and non-positively curved manifolds \cite{K}, hyperbolic groups \cite{CM}, 
	\item groups acting properly and isometrically on Hilbert spaces \cite{HK},
	\item groups acting properly and isometrically on bolic spaces \cite{KS}, 
	\item groups with finite asymptotic dimension \cite{Y1},
	\item groups coarsely embeddable into Hilbert spaces \cite{Y2}\cite{H}\cite{STY}, 
	\item groups coarsely embeddable into Banach spaces with property (H) \cite{KY}, 
	\item all linear groups and subgroups of all almost connected Lie groups \cite{GHW}, 
	\item  subgroups of the mapping class groups \cite{Ha}\cite{Ki}, 
	\item subgroups of $\operatorname{Out}(F_n)$, the outer automorphism groups of the free groups \cite{BGH},
	\item groups acting properly and isometrically on (possibly infinite dimensional) admissible  Hilbert-Hadamard spaces,
	in particular geometrically discrete subgroups of the group of volume preserving diffeomorphisms of any smooth compact manifold \cite{GWY}. 
\end{enumerate}
In the first three cases, an isometric action of a discrete group $\Gamma$ on a metric space $X$ is said to be \emph{proper} if for some $x\in X$,  $d(x, gx)\rightarrow \infty$ as $g \rightarrow \infty$,  i.e. for  any $x\in X$ and  any positive number $R>0$, there exists a finite subset $F$ of $\Gamma$ such that $ d(x, gx)> R$ if $g \in \Gamma - F$.

In a tour de force, Connes proved a striking theorem that the Novikov conjecture holds for higher signatures associated to Gelfand-Fuchs classes \cite{C1}. Connes, Gromov, and Moscovici proved the Novikov conjecture for higher signatures associated to Lipschitz group cohomology classes \cite{CGM}. Hanke-Schick and Mathai proved the Novikov conjecture  for higher signatures associated to group cohomology classes with degrees one and two \cite{HS}\cite{Ma}.

J. Rosenberg discovered an important application of the (rational) strong Novikov conjecture to the existence problem of Riemannian metrics with  positive scalar curvature \cite{R}. We refer to Rosenberg's survey \cite{R1} for recent development on this topic.

For the topological approach to the Novikov conjecture, we refer to the  articles
\cite{N, FH, CP, DFW,  FW1, W1}.

\subsection*{Acknowledgment}
The author  wish to thank Xiaoman Chen, Paul, Baum, Alain Connes, Misha Gromov, Guihua Gong,  Sherry Gong, Erik Guentner, Nigel Higson, Gennadi Kasparov, Vincent Lafforgue, Herv\'e Oyono-oyono, John Roe, Georges Skandalis,  Xiang Tang, Romain Tessera, Jean-Louis Tu,
Shmuel Weinberger, Rufus Willett,  Jianchao Wu, Zhizhang Xie for inspiring discussions on the Novikov conjecture. The author would like to thank Sherry Gong, Hao Guo, and  Slava Grigorchuk
for very helpful comments on this article.

\section{Non-positively curved groups and hyperbolic groups}
\label{ }

In this section, we give a survey on the work of A. Mishchenko, G. Kasparov,
A. Connes and H. Moscovici, G. Kasparov and G. Skandalis on the Novikov conjecture for non-positively curved groups and hyperbolic groups.

In  \cite{M}, A. Mishchenko introduced a theory of infinite dimensional Fredholm representations of discrete groups to prove the following theorem.

\begin{thm} The Novikov conjecture holds if the fundamental group of a manifold 
acts properly,  isometrically and cocompactly on a simply connected manifold with non-positive sectional curvature.
\end{thm}

In \cite{K}, G. Kasparov developed  a bivariant K-theory, called KK-theory, to prove the following theorem.

\begin{thm} The Novikov conjecture holds if the fundamental group of a manifold 
acts properly and isometrically  on a simply connected manifold with non-positive sectional curvature.
\end{thm}

In  the same article, G. Kasparov applied  the  above theorem to prove the following theorem.

\begin{thm} The Novikov conjecture holds if the fundamental group of a manifold 
is a discrete subgroup of a Lie group with finitely many connected components.
\end{thm}

A. Connes and H. Moscovici proved the following theorem using powerful techniques from noncommutative geometry \cite{CM}.

\begin{thm} The Novikov conjecture holds if the fundamental group of a manifold 
is a hyperbolic group.
\end{thm}

The theory of hyperbolic groups was developed by Gromov \cite{G}.
The proof of Theorem 2.5 uses Connes' theory of cyclic cohomology in a crucial way. Cyclic cohomology theory plays the role of de Rham theory in noncommutative geometry and  is the natural receptacle for the Connes-Chern character \cite{C}.

The following theorem of G. Kasparov and G. Skandalis unified the above results
\cite{KS}.

\begin{thm} The Novikov conjecture holds if the fundamental group of a manifold 
is bolic.
\end{thm}

Bolicity is a notion of non-positive curvature. Examples of bolic groups include hyperbolic groups and groups acting properly and isometrically on simply connected manifolds with non-positive sectional curvature.

\section{Amenable groups, groups with finite asymptotic dimension and coarsely embeddable groups} 
\label{}

In this section, we give a survey on the work of Higson-Kasparov 
on the Novikov conjecture for amenable groups, my work on the Novikov conjecture for groups with finite asymptotic dimension, and the work of G. Yu,
N. Higson, Skandalus-Tu-Yu on the Novikov conjecture for groups coarsely embeddable into Hilbert spaces. Finally we discuss the work of Kasparov-Yu on the connection of the Novikov conjecture with Banach space geometry.

Higson and Kasparov developed an index theory of certain differential operators on Hilbert space to prove the following theorem \cite{HK}.

\begin{thm} The Novikov conjecture holds if the fundamental group of a manifold 
acts properly  and isometrically on a Hilbert space.
\end{thm}

Recall that an isometric action $\alpha$ of a group $\Gamma$ on a Hilbert space is said to be proper if $ \|  \alpha(\gamma)h \| \rightarrow \infty$ when $\gamma \rightarrow \infty$, i.e. for  any $h\in H$ and  any positive number $R>0$, there exists a finite subset $F$ of $\Gamma$ such that $\| \alpha(\gamma) h \| > R$ if $\gamma \in \Gamma - F$.
A theorem of Bekka-Cherix-Valette states that every amenable group acts properly and isometrically on Hilbert space \cite{BCV}. Roughly speaking, a group is amenable if there exists large finite subsets of the group with small boundary.  The concept of amenability is a large-scale geometric property and  was introduced by von Neumann. We refer the readers to the book \cite{NY} as a general reference for geometric group theory related to the Novikov conjecture.

\begin{corr}
The Novikov conjecture holds if the fundamental group of a manifold is amenable.
\end{corr}

This corollary is quite striking since the geometry of amenable groups can be very complicated (for example, the Grigorchuk groups \cite{Gr}).

Next we recall a few basic concepts from geometric group theory.
A non-negative function $l$ on a countable group $G$ is called a length function if 
(1) $l(g^{-1})=l(g)$ for all $g\in G$; (2) $l(gh)\leq l(g) +l(h)$ for all $g$ and $h$ in $G$; (3) $l(g)=0$ if and only if $g=e$, the identity element of $G$.
We can associate a left-invariant length metric $d_l$ to $l$: $d_l(g,h)=l(g^{-1}h)$ for all $g,h\in G$. A length metric is called proper if the length function is a proper map (i.e. the inverse image of every compact set is finite in this case). 
It is not difficult to show that every countable group $G$ has a proper length metric. If $l$ and $l'$ are two proper length functions on $G$, then their associated length metrics are coarsely equivalent.  If $G$ is a  finitely generated group and $S$ is a finite symmetric generating set (in the sense that if an element is in $S$, then its inverse is also in $S$), then we can define the  word length $l_S$ on $G$ by
$$l_S(g)=min \{n: g=s_1\cdots s_n, s_i\in S\}.$$
If $S$ and $S'$ are two finite symmetric generating sets of $G$, then their associated proper length metrics are quasi-isometric.

The following concept is due to Gromov \cite{G1}.

\begin{defns} The asymptotic dimension of a proper metric space $X$ is the smallest integer $n$ such that for every $r>0$,
there exists a uniformly bounded cover $\{ U_i\}$ for which the number of $U_i$ intersecting each $r$-ball $B(x, r)$ is at most $n+1$.
\end{defns}

As example, the asymptotic dimension of ${\mathbb Z}^n$ is $n$, and the asymptotic dimension of the free group ${\mathbb F}_n$ with $n$ generators is $1$. Asymptotic dimension is invariant under coarse equivalence. The Lie group $GL(n, {\mathbb R})$ with a left-invariant Riemannian metric is quasi-isometric to $T(n, {\mathbb R})$, the subgroup of invertible upper triangular matrices.
By permanence properties of asymptotic dimension [BD1], we know that the solvable group $T(n, {\mathbb R})$ has finite asymptotic dimension. As a consequence, every countable discrete subgroup of $GL(n, {\mathbb R})$ has finite asymptotic dimension
(as a metric space with a proper length metric).
More generally, one can prove that every discrete subgroup of an almost connected Lie group has finite asymptotic dimension 
(a Lie group is said to be almost connected if the number of its connected components is finite). 
Gromov's hyperbolic groups  have finite asymptotic dimension \cite{Roe2}.
Mapping class groups also have finite asymptotic dimension \cite{BBF}.

In \cite{Y1}, I developed a quantitative  operator K-theory to prove the following theorem.

\begin{thm} The Novikov conjecture holds if the fundamental group of a manifold 
has finite asymptotic dimension.
\end{thm}

The following concept  of Gromov makes precise the idea of drawing a good picture of a metric space in a Hilbert space.

\begin{defns}(Gromov): Let $X$ be a metric space and  $H$ be a Hilbert space.  A map $f : X \rightarrow
H$ is said to be a coarse embedding if there exist non-decreasing
functions $\rho_1$ and $\rho_2$ on $[0,  \infty )$ such that
\item {(1)} $\rho_1 (d(x,y)) \leq d_H (f(x), f(y)) \leq \rho_2 (d(x,y))$
for all $x,y \in X$;
\item {(2)} $\lim_{r \rightarrow +\infty} \rho_1 (r) = + \infty$.

\end{defns}

Coarse embeddability of a countable group is independent of the choice of proper length metrics. 
Examples of groups coarsely embeddable into Hilbert space include groups acting properly and isometrically on a Hilbert space (in particular amenable groups \cite{BCV}), groups with Property A \cite{Y2}, 
 countable subgroups of connected Lie groups \cite{GHW}, 
 hyperbolic groups \cite{S},  groups with finite asymptotic dimension, Coxeter groups \cite{DJ}, mapping class groups \cite{Ki, Ha}, 
 and semi-direct products of groups of the above types.

The following theorem unifies the above theorems.

\begin{thm} The Novikov conjecture holds if the fundamental group of a manifold 
is coarsely embeddable into Hilbert space.
\end{thm}

Roughly speaking, this theorem says if we can draw a good picture of the fundamental group in a Hilbert space, then we can recognize the manifold at an infinitesimal level. This theorem was proved by myself when the classifying space of the fundamental group has the homotopy type of a finite CW complex \cite{Y2} and this finiteness condition was removed by N. Higson \cite{H}, Skandalis-Tu-Yu \cite{STY}. The original proof of the above result makes heavy use of infinite diimensional analysis.  More recently, R. Willett and myself found a relatively elementary proof within the framework of basic operator K-theory
\cite{WiY}.

E. Guentner, N. Higson and S. Weinberger proved the beautiful theorem that
linear groups are coarsely embeddable into Hilbert space \cite{GHW}.  Recall that a group is called linear if it is a subgroup of $GL(n, k)$ for some field $k$. The following theorem follows as a consequence \cite{GHW}.

\begin{thm} The Novikov conjecture holds  if the fundamental group of a manifold is a linear group.
\end{thm}

More recently, Bestvina-Guirardel-Horbez proved that 
$\operatorname{Out}(F_n)$, the outer automorphism groups of the free group $F_n$,
is coarsely embeddable into Hilbert space. This implies the following theorem \cite{BGH}.

\begin{thm} The Novikov conjecture holds   if the fundamental group of a manifold is a subgroup of  $\operatorname{Out}(F_n)$.
\end{thm}

We have the following open question.

\begin{open} Is  every countable subgroup of the diffeomorphism group of the circle coarsely embeddable into Hilbert space?
\end{open}

Let $\mathfrak{E}$ be the smallest class of groups which include all groups coarsely embeddable into Hilbert space and is closed under direct limit.  Recall that if $I$ is a directed set and $\{G_i, \phi_{i,j}\}_{i,j\in I, i< j}$ is a direct system of groups
over $I$, then we can define the direct limit $\lim G_i$, where 
the homomorphism $\phi_{ij}: G_i\rightarrow G_j$ for $i< j$ is not necessary injective. 

The following result is a consequence of Theorem 3.3.

\begin{thm} The Novikov conjecture holds   if the fundamental group of a manifold is in the class $\mathfrak{E}$.
\end{thm}

The following open question is a challenge to geometric group theorists.

\begin{open} Is there any countable group not in the class $\mathfrak{E}$?
\end{open}

We mention that the Gromov monster groups are in class $\mathfrak{E}$ \cite{G2, G3, AD, O}.

Next we shall discuss the connection of the Novikov conjecture with the geometry of Banach spaces.

\begin{defns}  A  Banach space $X$ is said to have Property (H) if there exist  an increasing sequence of finite dimensional subspaces
$\{ V_n\} $ of $X$ and an increasing  sequence of finite dimensional subspaces $\{ W_n\} $ of a Hilbert space such that
\begin{enumerate}
\item[(1)]
 $V=\cup_n V_n$ is dense in $X$,

\item[(2)]  if  $W=\cup_n W_n$,  $S(V)$ and $S(W)$ are respectively the unit spheres of $V$ and $W$, then there exists a uniformly continuous map $\psi: S(V)\rightarrow S(W)$ such that
the restriction of $\psi$ to $S(V_n)$ is a homeomorphism (or more generally a degree one map) onto $S(W_n)$ for each $n.$

\end{enumerate}

\end{defns}

As an example, let $X$ be the  Banach space $l^p(\mathbb{N})$ for some  $p\geq 1$.
Let $V_n$  and $W_n$  be respectively the subspaces of $l^p(\mathbb{N})$  and $l^2(\mathbb{N})$ consisting of all sequences whose coordinates are zero after the  $n$-th terms.
We define a map $\psi$ from $S(V)$ to $S(W)$ by
 $$\psi(c_1, \cdots, c_k, \cdots)= (c_1 |c_1|^{p/2-1}, \cdots, c_k |c_k|^{p/2-1},\cdots). $$
 $\psi$ is called the Mazur map.
It is not difficult to verify that $\{V_n\}, \{W_n\}$ and $\psi$ satisfy the conditions in the definition of Property (H).
For each $ p\geq 1$, we can similarly prove that $C_p$, the Banach space of all Schatten $p$-class operators on a Hilbert space, has Property (H).  

Kasparov and Yu proved the following.

\begin{thm} The Novikov conjecture holds if the fundamental group of a manifold 
is coarsely embeddable into a Banach space with Property (H).
\end{thm}

Let $c_0$ be the Banach space consisting of  all sequences of real numbers converging  to $0$ with the sup norm.

\begin{open} Does  the Banach space $c_0$ have Property (H)?
\end{open}

A positive answer to this question would imply the Novikov conjecture since every countable group admits a coarse embedding into $c_0$ \cite{BG}.

A less ambitious question is the following.

\begin{open}
Is every countable subgroup of the diffeomorphism group of a compact smooth manifold  coarsely embeddable into $C_p$ for some $p\geq 1$?
\end{open}

For each $p>q\geq 2$, it is also an open question to construct a bounded geometry space
which is coarsely embeddable into $l^p(\mathbb{N})$ but not $l^q(\mathbb{N})$.  Beautiful  results in \cite{JR} and \cite{MN}
indicate that  such a construction should be possible. Once such a metric space is constructed, the next natural question is to construct countable groups which coarsely contain such a metric space. These groups would be from another universe and would be different from any group we currently know.

\section{Gelfand-Fuchs classes, the group of volume preserving diffeomorphisms, Hilbert-Hadamard spaces}
\label{}

In this section, we give an overview on the work of A. Connes, Connes-Gromov-Moscovici on the Novikov conjecture for Gelfand-Fuchs classes and the recent work of  Gong-Wu-Yu on the Novikov conjecture for groups acting  properly and isometrically on  Hilbert-Hadamard spaces and for any geometrically discrete subgroup of the group of volume preserving diffeomorphisms of a compact smooth manifold.

A. Connes proved the following deep theorem on the Novikov conjecture \cite{C1}.

\begin{thm} 
The Novikov conjecture holds for higher signatures associated to the Gelfand-Fuchs cohomology classes of a subgroup of the group of diffeomorphisms of a compact smooth manifold.
\end{thm}

The proof of this theorem uses the full power of noncommutative geometry \cite{C}.

\begin{open} Does the Novikov conjecture hold for any subgroup of the group of diffeomorphisms of a compact smooth manifold?
\end{open}

Motivated in part by this open question, S. Gong, J. Wu and G. Yu prove the following theorem \cite{GWY}.

\begin{thm} 
The Novikov conjecture holds for groups acting properly and isometrically on an admissible Hilbert-Hadamard space. 
\end{thm}

Roughly speaking,  Hilbert-Hadamard spaces are (possibly infinite dimensional) simply connected spaces with non-positive curvature. We will give a precise definition a little later.
We say that a  Hilbert-Hadamard  space $M$ is \textit{admissible} if it has a sequence of subspaces $M_n$, whose union is dense in $M$, such that each $M_n$, seen with its inherited metric from $M$, is isometric to a finite-dimensional Riemannian manifold. Examples of admissible Hilbert-Hadamard spaces
include all simply connected and non-positively curved Riemannian manifold, the Hilbert space, and certain infinite dimensional symmetric spaces. Theorem 4.3
can be viewed as a generalization of both Theorem 2.1 and Theorem 3.1.

Infinite dimensional symmetric spaces are often naturally admissible Hilbert-Hadamard spaces. One such an example 
 of an admissible infinite-dimensional symmetric space is 
\[L^2(N,\omega, \operatorname{SL}(n, \mathbb{R})/\operatorname{SO}(n)),\]
where $N$ is a compact smooth manifold with a given volume form $\omega$.
This infinite-dimensional symmetric space
 is defined to be the completion of  the space of all smooth maps 
from $N$ to $X=\operatorname{SL}(n, \mathbb{R})/\operatorname{SO}(n)$
with respect to the following distance:
		$$d(\xi, \eta) = \left( \int_N (d_X(\xi(y),\eta(y)))^2 \, d\omega(y) \right)^{\frac{1}{2}},$$ where $d_X$ is the standard Riemannian metric on the symmetric space $X$ and $\xi$ and $\eta$ are two smooth maps from $N$ to $X$. This space can be considered as  the space of $L^2$-metrics on $N$ with the given volume form 
$\omega$  and is Hilbert-Hadamard space.
With the help of this infinite dimensional symmetric space, the above theorem can be applied to study the Novikov conjecture for geometrically discrete subgroups of the group of volume preserving diffeomorphisms on such a manifold.

The key ingredients of the proof for Theorem 4.3 include  a construction of a $C^\ast$-algebra modeled after the Hilbert-Hadamard  space, a deformation technique for the isometry group of the Hilbert-Hadamard  space and its corresponding actions on $K$-theory, and a $KK$-theory with real coefficient developed by Antonini, Azzali, and Skandalis \cite{AAS}. 

Let $\operatorname{Diff}(N,\omega)$ denote the group of volume preserving diffeomorphisms on a compact orientable smooth manifold $N$ with a given volume form $\omega$. In order to define the concept of geometrically discrete subgroups of $\operatorname{Diff}(N,\omega)$, let us fix a Riemannian metric on $N$  with the given volume  $\omega$ and define a length function $\lambda$ on  $\operatorname{Diff}(N,\omega)$ by:
\[ \lambda_+(\varphi) = \left(\int_N (\log(\|D\varphi\|))^2 d\omega\right)^{1/2}\]
and
\[ \lambda(\varphi) =  \max\left\{\lambda_+(\varphi),\lambda_+(\varphi^{-1})\right\} \; \]
 for all $\varphi \in \operatorname{Diff}(N,\omega)$, 
where $D\varphi$ is the Jacobian of $\varphi$, and the norm $\|\cdot \|$ denotes the operator norm, computed using the chosen Riemannian metric on $N$.

\begin{defns}
A subgroup $\Gamma$ of $\operatorname{Diff}(N,\omega)$ is said to be a geometrically discrete subgroup 
if $\lambda(\gamma) \rightarrow \infty$ when $\gamma \to \infty$ in a $\Gamma$, i.e. for any $R>0$, there exists a finite subset $F\subset \Gamma$ such that $\lambda(\gamma)\geq R$ if $\gamma \in \Gamma \setminus F$. 
\end{defns}

Observe that although the length function $\lambda$ depends on our choice of the Riemannian metric, the above notion does not. Also notice that if $\gamma$ preserves the Riemannian metric we chose, then $\lambda(g)=0$. This suggests that the class of geometrically discrete subgroups of $\operatorname{Diff}(N,\omega)$  doesn't intersect with the class of groups of isometries. Of course we already know the Novikov conjecture for any group of isometries on a compact Riemannian  manifold. This, together with the following result, gives an optimistic perspective on the open question on the Novikov conjecture for groups of volume preserving diffeomorphisms.

\begin{thm} 
	Let $N$ be a compact smooth manifold with a given volume form $\omega$, and let $\operatorname{Diff}(N,\omega)$ be the group of all volume preserving diffeomorphisms of $N$. The  Novikov conjecture holds for any geometrically discrete subgroup of $\operatorname{Diff}(N,\omega)$. 

\end{thm}

Next we will give a precise definition of Hilbert-Hadamard space. 

We will first recall the concept of CAT(0) spaces. 
Let $X$ be a geodesic metric space.
Let $\Delta $ be a triangle in $X$ with geodesic segments as its sides. $\Delta$
  is said to satisfy the CAT(0) inequality if there is a comparison triangle 
  $\Delta '$ in Euclidean space, with sides of the same length as the sides of 
  $ \Delta$ , such that distances between points on $\Delta $ are less than or equal to the distances between corresponding points on $ \Delta '.$
  The geodesic metric space $X$ is said to be a CAT(0) space if every geodesic triangle satisfies the CAT(0) inequality. 

Let $X$ be a geodesic metric space. For three distinct points $x,y,z \in X$, we define the comparison angle $\widetilde{\angle} xyz$ to be
\[\widetilde{\angle} xyz = \arccos\left(\frac{d(x,y)^2 + d(y,z)^2 - d(x,z)^2}{2d(x,y)d(y,z)}\right).\]
In other words, $\widetilde{\angle} xyz$ can be thought of as the angle at $y$ of a triangle $ABC$ in the Euclidean plane with side-lengths that agree with the side-lengths of the geodesic triangle $xyz$ in $X$.

Given two nontrivial geodesic paths $\alpha$ and $\beta$ emanating from a point $p$ in $X$, meaning that $\alpha(0) = \beta(0) = p$, we define the angle between them, $\angle(\alpha,\beta)$, to be 
\[\angle(\alpha,\beta) = \lim_{s,t \to 0} \widetilde{\angle} (\alpha(s),p,\beta(t)) \; ,\]
provided that the limit exists. For CAT(0) spaces, since the comparison angle $\widetilde{\angle} (\alpha(s),p,\beta(t))$ decreases with $s$ and $t$, the angle between any two geodesic paths emanating from a point is well-defined. The above angles satisfy the triangle inequality.

For a point $p \in X$, let $\Sigma_p'$ denote the metric space induced from the space of all geodesics emanating from $p$ equipped with the pseudometric of angles, that is, for geodesics $\alpha$ and $\beta$, we define $d(\alpha,\beta) =\angle(\alpha,\beta)$. Note, in particular, from our definition of angles, that $d(\alpha,\beta) \leq \pi$ for any geodesics $\alpha$ and $\beta$.

We define  $\Sigma_p$ to be the completion of $\Sigma_p'$ with respect to the distance $d$. 
The \textit{tangent cone} $K_p$ at a point $p$ in $X$ is then defined to be a metric space which is, as a topological space, the cone of $\Sigma_p$. That is, topologically
\[K_p \simeq \Sigma_p \times [0,\infty)/\Sigma_p \times \{0\}.\]

The metric on it is given as follows. Any  two points $p,q \in K_p$ can expressed   as $p=[(x,t)]$ and $q = [(y,s)]$. The metric is given by
\[d(p,q) = \sqrt{t^2+s^2-2st\cos(d(x,y))}.\]
The distance is what the distance would be if we went along geodesics in a Euclidean plane with the same angle between them as the angle between the corresponding directions in $X$.

The following definition is inspired by \cite{FS}.

\begin{defn}\label{defn:hhs}
	A \emph{\hhs} is a complete geodesic CAT(0) metric space (i.e., an Hadamard space) all of whose tangent cones are isometrically embedded in Hilbert spaces. 
\end{defn}

	Every connected and simply connected \emph{Riemannian-Hilbertian manifold with non-positive sectional curvature} is a separable {\hhs}. In fact, a Riemannian manifold without boundary is a {\hhs} if and only if it is complete, connected, and simply connected, and has nonpositive curvature. 
We remark that a CAT(0) space $X$ is always uniquely geodesic. 

Recall that a subset of a geodesic metric space is called \emph{convex} if it is again a geodesic metric space when equipped with the restricted metric. We observe that a closed convex subset of a {\hhs} is itself a {\hhs}.

\begin{defns}
	A separable {\hhs} $M$ is called \emph{admissible} if there is a sequence of convex subsets isometric to finite-dimensional Riemannian manifolds whose union is dense in $M$. 
\end{defns}

The notion of {\hhs}s is more general than simply connected  Riemannian-Hilbertian space with non-positive sectional curvature. For example, the infinite dimensional symmetric space $L^2(N,\omega, \operatorname{SL}(n, \mathbb{R})/\operatorname{SO}(n))$
is a {\hhs} but not a Riemannian-Hilbertian space with non-positive sectional curvature.
	
\section{Geometric complexity, topological rigidity and non-rigidity}

An integral version of the Novikov conjecture implies the stable Borel conjecture \cite{J, FH} which states that all compact aspherical manifolds
are stably topologically rigid, i.e. if another compact manifold $N$ is homotopy equivalent to the given compact aspherical manifold $M$, then $N\times \mathbb{R } ^n$ is homeomorphic to $M\times \mathbb{R}^n$ for some $n$. 
In this section, we discuss the concept of decomposition complexity introduced by Guentner-Tessera-Yu [GTY] and  its applications to the stable Borel conjecture. The method used here is very much inspired by the quantitative operator K-theory approach developed in  \cite{Y1}. We shall also briefly discuss the application of the Novikov conjecture to non-rigidity of manifolds in Weinberger-Xie-Yu  
 \cite{WXY}.

We shall first recall the concept of finite decomposition complexity.

 For any $r>0$, a collection of subspaces $\{\, Z_i \,\}$ of a metric
space $Z$ is said to be $r$-disjoint if for all $i\neq j$ we have
$d(Z_i,Z_j)\geq r$. To express the idea that $Z$ is the union of
subspaces $Z_i$, and that the collection of these subspaces is
$r$-disjoint we write \begin{equation*}
  Z = \bigsqcup_{r-{\rm disjoint}} Z_i.
\end{equation*}
A family of  metric spaces $\{\, Z_i \,\}$ is
called  bounded if there is a uniform bound on the diameter of the
individual $Z_i$:
\begin{equation*}
  \sup {\rm diameter}(Z_i) <\infty.
\end{equation*}

\begin{defns}
A  family  of metric spaces $\{ X\} $ is  $r$-decomposable over another family of  metric
spaces  $\{ Y\} $ if every $X\in \{ X\} $ admits an {\it $r$-decomposition\/}
\begin{equation*}
    X = X_0\cup X_1, \quad
     X_i = \bigsqcup_{r-{\rm disjoint}} X_{ij},
\end{equation*}
where each $X_{ij}\in \{ Y\}$.  

\end{defns}

\begin{defns}
  Let $\Omega $ be a collection of families of metric spaces.  A  family of metric spaces $\{X\}$ is said to be decomposable over $\Omega$ if, for every $r>0$,
 there exists a family of   metric spaces $\{ Y \} \in\Omega$ and an $r$-decomposition of
 $\{ X\} $ over $\{Y\}$.
The collection $\Omega$ is said to be stable under decomposition if every
  family  of metric spaces which decomposes over $\Omega$ actually belongs to $\Omega$.

\end{defns}

\begin{defns} 
  The collection $\mathfrak{D}$ of  families of metric spaces with finite
  decomposition complexity is the minimal collection of families of metric spaces
  containing all bounded families of metric spaces and stable under decomposition.
  We abbreviate membership in $\mathfrak{D}$ by saying that a  family of metric spaces in
  $\mathfrak{D}$ has finite decomposition complexity. A metric space $X$ is said to have finite decomposition complexity
  if the family consisting of only $X$ has finite decomposition complexity.
\end{defns}

Observe that finite decomposition complexity is invariant under coarse equivalence. We also remark that metric spaces with finite decomposition complexity are coarsely embeddable into Hilbert space. 
 
 By the definition of asymptotic dimension, any proper metric space with asymptotic dimension at most $1$ has finite decomposition complexity. 
 More generally,  a metric space with  finite asymptotic dimension has finite decomposition complexity.
 This fact follows from a theorem of Dranishnikov and Zarichnyi stating that a proper metric space with finite asymptotic dimension is coarsely equivalent to a subspace of  the product of finitely many trees.
 
 Let $G=\oplus_{k=1}^{\infty} {\mathbb Z}$ be the countable group with the proper length metric associated to the length 
 function $l$: 
 $$l(\oplus_{k=1}^{\infty} n_k)= \sum_{k=1}^{\infty} k |n_k| $$ for each $\oplus_{k=1}^{\infty} n_k\in G$.
 For each $r>0$, let $k_0$ be the smallest integer greater  than $r$.
 For each $\alpha \in \oplus_{k=k_0}^{\infty} {\mathbb Z}$, let 
 $$G_\alpha=\{ \oplus_{k=1}^{k_0-1} n_k \oplus  \alpha: \oplus_{k=1}^{k_0-1} n_k \in \oplus_{k=1}^{k_0-1} {\mathbb Z}\}.$$
 Notice that $G=\bigsqcup_{r-{\rm disjoint}} G_\alpha$ and that $G_\alpha$  is (uniformly) coarsely equivalent to 
 $\oplus_{k=1}^{k_0-1} {\mathbb Z}$.  This implies  that  $\{ G_\alpha\}\in \mathfrak{D}$.
 It follows that $G$ has finite decomposition complexity despite the fact that $G$ has infinite asymptotic dimension.
 If $H \subset GL(2, {\mathbb R})$ is  the finitely generated  group consisting of all matrices of the form
 $\left(\begin{array}{cc} \pi^n&p(\pi) \\ 0&\pi^{-n} \end{array}\right)$, with a little extra work we can show that $H$ has finite decomposition complexity.
 
More generally, we have the following result from [GTY].

\begin{thm}  Any countable subgroup of $GL(n,k)$  has finite decomposition complexity (as a metric space with a proper length metric), where $k$ is a field.
\end{thm}

The same result is true for any countable subgroup of an almost connected Lie group  and any countable elementary amenable group [GTY].

The following result is proved in [GTY].

\begin{thm}
The stable Borel conjecture holds for aspherical manifolds whose fundamental groups have finite decomposition complexity.
\end{thm}

We have the following open question.

\begin{open}
Does any  countable amenable group have finite decomposition complexity?
\end{open}

In particular, it remains an open question whether the Grigorchuk groups
have finite decomposition complexity. 

Finally, we ask the following open question.

\begin{open}
Does the group $\operatorname{Out}(F_n)$  have finite decomposition complexity?
\end{open}

There has been  spectacular recent progress on the Borel conjecture. We will not attempt to survey all important results, but only mention
the fundamental work of Farrell-Jones \cite{FJ1, FJ2, FJ3, FJ4}, Bartels-L\"{u}ck \cite{BL},  and Bartels-L\"{u}ck-Holger \cite{BLH}.  I also refer interested readers to the beautiful books by T. Farrell and L. Jones \cite{FJ2} and  S. Weinbeger \cite{W}.

Finally we mention interesting applications of the Novikov conjecture to non-rigidity of manifolds \cite{WY, WXY}. In particular, Weinberger-Xie-Yu introduced the Novikov rho invariant to prove that the structure groups of certain manifolds are infinitely generated \cite{WXY}. Recall that the topological structure group $S^{TOP}(M)$ is the abelian group of equivalence classes of all
pairs $(f, X) $ such that $X$ is a closed oriented manifold and $f :  X\rightarrow M$ is an
orientation-preserving homotopy equivalence. The structure group $S^{TOP}(M)$
measures non-rigidity of $M$.

\begin{thm} Let $M$ be a closed oriented topological manifold of dimension $n \geq 5$, and $\Gamma$ be its fundamental group. Suppose the rational Kasparov-Baum-Connes assembly map for $\Gamma$ is
split injective. If $\oplus_{k\in \mathbb{Z}} H_{n+1+4k} (\Gamma, \mathbb{C})$
is infinitely generated, then the topological structure group of $S^{TOP}(M)$
is infinitely generated.
\end{thm}

We refer to the article \cite{WXY} for examples of groups satisfying conditions in the above theorem. We remark that the condition on split injectivity of the rational Kasparov-Baum-Connes assembly map  is one version of the strong Novikov conjecture.

\bigskip
\footnotesize

\end{document}